# LA ESPIRAL AUREA, SU LONGITUD Y RECTANGULOS AUREOS

CAMPO ELIAS GONZALEZ PINEDA.

cegp@utp.edu.co

SANDRA MILENA GARCIA.

tazyotas@utp.edu.co

## RESUMEN

En este artículo estudiamos la espiral aurea y calculamos su longitud. Estudiamos los rectángulos áureos y calculamos algunas medidas de interés. Veremos como los únicos rectángulos que pueden subdividirse en rectángulos indefinidamente son los áureos. Además, mostraremos que algunos rectángulos se asemejan con un cierto grado k en ser áureos.

## Rectángulos áureos y algunos cálculos de interés

En este artículo estudiamos algunas relaciones que se obtienen al dividir un rectángulo en subrectángulos más pequeños. Sabemos que un rectángulo de lados $a > b$ es áureo si

$$\frac{a}{b} = \varphi = \frac{1+\sqrt{5}}{2}$$

Sea $R$ el rectángulo, si a éste le extraemos el cuadrado de lado $b$ se obtiene el rectángulo $R_1$ de lados $b, a-b$. Es claro que este rectángulo también es áureo. En efecto,

$$\frac{b}{a-b} = \frac{b/a}{1-b/a} = \frac{1/\varphi}{1-1/\varphi} = \varphi$$

Recuérdese que

$$\varphi^2 = \varphi + 1 \iff \varphi = 1 + \frac{1}{\varphi} \iff \varphi - 1 = \frac{1}{\varphi}$$

Si continuamos el proceso se obtiene una sucesión infinita de rectángulos áureos. Nuestro interés es tratar de encontrar ciertos números asociados a estos rectángulos y sus similares.

Se sabe que en el rectángulo áureo aparece una espiral que resulta de unir los arcos de circunferencia del cuadrado mayor dibujado dentro del rectángulo.

En este artículo mediremos la longitud de la espiral y algunas áreas relacionadas dentro de dichos rectángulos. Observemos que si

$$\frac{a}{b} = \varphi \iff a = \varphi b$$

Como veremos, al subdividir cada rectángulo en rectángulos más pequeños aparece de manera natural la sucesión de Fibonacci.

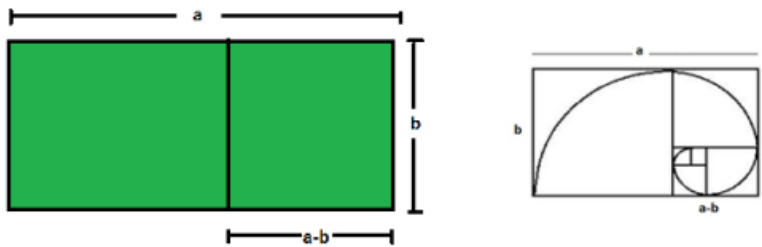

**Nota**. En todo este artículo supondremos rectángulos de lados desiguales, es decir, no cuadrados. Diremos que un rectángulo es horizontal si su lado más grande es horizontal y vertical si su lado más largo es vertical. Así un rectángulo áureo se puede dividir en un cuadrado y un rectángulo áureo vertical al dibujarle el cuadrado más grande.

Para el caso del rectángulo áureo decimos que éste se puede subdivir o que genera infinitos rectángulos áureos. Como veremos esta es una cualidad especial de estos rectángulos, en el sentido de que son los únicos que tienen esta propiedad. Sin embargo, algunos rectángulos no pueden subdividirse en rectángulos. Para el caso de rectángulos de lados $a = 2l, \ b = l$, es decir, si $\frac{a}{b} = 2$ se divide en dos cuadrados de lado $l$. Otros tienen la propiedad de que pueden subdivirse en un número finito de subrectángulos. Es de resaltar que los rectángulos aureos generan una sucesion intercalada de rectángulos horizontales y verticales, que al dibujarlos en su interior genera la famosa espiral aurea.

## Regiones de interés en un rectángulo

Llamaremos $L_k$ el arco de circunferencia, $A_k$ el área del cuarto de circunferencia (región amarilla) correspondiente, $B_k$ el área del rectángulo menor que resulta (región morada), $C_k$ el área de la región azul, $D_k$ la medida de la diagonal del cuadrado mayor del rectángulo, ver figura.

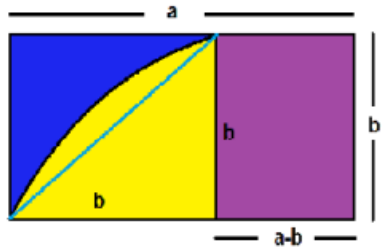

Para este propósito es necesario recodar la sucesión de Fibonacci algunas de sus propiedades y también el número $\varphi$ del cual también recordaremos algunas de sus propiedades.

## Sucesión de Fibonacci.

La sucesión de Fibonacci se define por la fórmula de recurrencia,

$$F_1 = 1, F_2 = 1, \ y \ \ F_n = F_{n-1} + F_{n-2} \text{ Para } n \geq 3.$$

Enlistados algunos términos tenemos.

$$F_1 = 1, F_2 = 1, F_3 = 2, F_4 = 3, F_5 = 5, F_6 = 8, F_7 = 13, F_8 = 21,$$

## Algunas propiedades

Algunas propiedades de la sucesión de Fibonacci se enuncian sin demostración, éstas se pueden consultar en la bibliografía

1. Dos números consecutivos de Fibonacci siempre son primos relativos.
2. Para cualesquiera números de Fibonacci consecutivos $A, B, C, D$ se tiene que $C^2 = AD + B^2$.
3. A excepción del 3 todo número de Fibonacci que sea primo tiene subíndice primo.
4. En el triángulo de Tartaglia (Pascal) sumando los términos de las diagonales secundarias, obtenemos los términos de la sucesión de Fibonacci.
5. $F_1 + F_2 + \cdots + F_n = F_{n+2} - 1$
6. $F_1 + F_3 + F_5 + \cdots + F_{2n-1} = F_{2n}$
7. $F_2 + F_4 + \cdots + F_{2n} = F_{2n+1} - 1$
8. $F_1 - F_2 + F_3 - F_4 + \cdots + (-1)^{n+1}F_n = (-1)^{n+1}F_{n-1} + 1$
9. $F_1 - F_2 + F_3 - F_4 + \cdots + F_{2n-1} - F_{2n} = -F_{2n-1} + 1$
10. $F_1 - F_2 + F_3 + \cdots + F_{2n-1} - F_{2n} + F_{2n+1} = F_{2n} + 1$
11. $F_1^2 + F_2^2 + \cdots + F_n^2 = F_n \cdot F_{n+1}$
12. $F_1 \cdot F_2 + F_2 \cdot F_3 + \cdots + F_{2n-1} \cdot F_n = F_{2n}^2$
13. $F_1 \cdot F_2 + F_2 \cdot F_3 + \cdots + F_{2n} \cdot F_{2n+1} = F_{2n+1}^2 - 1$
14. $F_1 + 2F_2 + \cdots + nF_n = nF_{n+2} - F_{n+2} + 2$
15. $F_n^2 + F_{n+1}^2 = F_{2n+1}$
16. $F_{n-1} \cdot F_{n+1} = F_n^2 + (-1)^n$
17. $F_{n+m} = F_{n-1}F_m + F_nF_{m+1}$
18. $F_{n+1}^2 - F_{n-1}^2 = F_{2n}$

Existen otras muchas propiedades que se dejan a la curiosidad del lector, pero quizá es conveniente resaltar que

$$\lim_{n \to \infty} \frac{F_{n+1}}{F_n} = \varphi$$

Donde como sabemos

$$\varphi = \frac{1 + \sqrt{5}}{2} = 1,618033...$$

**El número áureo o divina proporción**

El número áureo divina proporción tiene un largo historial, ameno y curioso que puede ser consultado en la bibliografía. Aquí nos limitamos a un tratamiento algebraico y damos algunas propiedades básicas, que nos sirven para ilustrar los procesos que veremos más adelante.

Fi, conocido como la constante relacionada con la belleza aun siendo irracional y nada bello como pensaba Pitágoras.

Utilizado al parecer desde el año 2000 ac. Lastimosamente como un supuesto ya que no existe evidencia que lo pruebe, los antiguos al parecer lo incluían en sus cálculos para las construcciones.

A este maravilloso número se le atribuyen desde siempre propiedades mágicas y místicas y, es tal vez por esa misma razón que no ocupa el mismo lugar dentro del cálculo que otros números irracionales como $\pi, e$ o $\sqrt{2}$. Fi es el único número real que satisface la ecuación,

$$\varphi^2 = \varphi + 1, \quad \varphi - 1 = \frac{1}{\varphi}, \quad \varphi^3 = \frac{\varphi + 1}{\varphi - 1}$$

Las potencias del número áureo pueden ser escritas en función de una suma de potencias de grados inferiores del mismo número, estableciendo una verdadera sucesión recurrente de potencias.

$$\varphi^n = \sum_{i=0}^{k/2} \binom{k/2}{i} \varphi^{n-(\frac{k}{2}+i)}$$

donde $i \in N_k$, k un número par.

El número áureo $\varphi$ es la unidad fundamental «e» del cuerpo $Q(\sqrt{5})$ con recíproco $-\bar{\varphi}$

En esta extensión el «emblemático» número irracional $\sqrt{2}$ cumple las siguientes igualdades:

$$\sqrt{2} = \varphi\sqrt{2(1+\bar{\varphi})} = -\bar{\varphi}\sqrt{2(1+\varphi)}$$

En realidad notemos que $\varphi^2 - \varphi = 1$ o también, $\varphi^2 + \varphi^2\bar{\varphi} = 1$ por lo que para todo real $p$ se tiene que

$$p = \varphi^2 p(1+\bar{\varphi})$$

El numero fi tiene historia en muchos campos de estudio, el arte, las matemáticas, la música, arquitectura, historia y en la misma naturaleza.

El primero en hacer un estudio formal sobre el número áureo fue Euclides, unos tres siglos antes de Cristo, en su obra Los Elementos. Casi 2000 años más tarde, en 1525, Alberto Durero publicó su Instrucción sobre la medida con regla y compás de figuras planas y sólidas, años más tarde también Johannes Kepler desarrolló su modelo del Sistema Solar, explicado en su libro El Misterio Cósmico, donde le da importancia a este número con su frase "La geometría tiene dos grandes tesoros: uno es el teorema de Pitágoras; el otro, la división de una línea entre el extremo y su proporcional". Además de la estrecha relación que guarda el número áureo con la famosa serie de Fibonacci. Si llamamos $F_n$ al enésimo número de Fibonacci y $F_{n+1}$ al siguiente, podemos ver que a medida que n se hace más grande, la razón entre

$F_{n+1}$ y $F_n$ oscila, siendo alternativamente menor y mayor que la razón áurea.

En el campo del arte, es tal vez Leonardo Da Vinci quien es más conocido por hacer uso de la proporción áurea en su obra maestra "**la monalisa**' y "**el hombre de Vitrubio**" por hablar de algunas. También lo hace su más cercano competidor, Miguel Ángel en "**el David"**

En la arquitectura, se dice que tanto en la antigüedad como en tiempos actuales la proporción áurea está presente en las construcciones del hombre, desde el Partenón, las pirámides de Egipto, la torre Eiffel hasta el desarrollo de un nuevo sistema de medida conocido como **"el modulor".**

Para terminar haremos una referencia de la proporción áurea presente en la naturaleza, comencemos con el ADN que es una código con toda la información acerca de un individuo, su conformación es espiral. Cualquiera de sus hélices, de 34 amstrongs de longitud y 21 de anchura, conforma un rectángulo áureo. Estos números son miembros de la sucesión de Fibonacci y el cociente entre ellos es aproximadamente $\varphi$.

En los vegetales como brócolis y coliflores se dibuja la espiral áurea y el plantas como los girasoles sus hoja están dispuestas de tal forma que cada parte este en contacto con la luz y los nutrientes necesarios.

## Para finalizar acerca de phi

Sabemos que $\varphi = \frac{1+\sqrt{5}}{2}, \ \bar{\varphi} = \frac{1-\sqrt{5}}{2},$ entonces,

1. $\varphi + \bar{\varphi} = 1$
2. $\varphi \cdot \bar{\varphi} = -1$
3. $\varphi^2 + \bar{\varphi}^2 = 3$
4. $\varphi^2 + \varphi\bar{\varphi} = \varphi$

También tenemos que

1. $\varphi^2 = \varphi + 1$
2. $\varphi^3 = 2\varphi + 1$
3. $\varphi^4 = 3\varphi + 2$
4. $\varphi^5 = 5\varphi + 3$
5. $\varphi^6 = 8\varphi + 3$
6. Por inducción $\varphi^n = F_n\varphi + F_{n-1}$ para $n \geq 2$

Sabiendo que $\varphi^n = F_n\varphi + F_{n-1}$ tenemos que $\bar{\varphi}^n = \bar{\varphi}F_n + F_{n-1}$. Sumando y restando encontramos

$$\varphi^n + \bar{\varphi}^n = F_{n+1} + F_{n-1}$$

$$\varphi^n - \bar{\varphi}^n = \sqrt{5}F_n$$

# Rectángulos y su relación con la sucesión de Fibonacci y el número áureo.

En esta parte realizaremos la división de un rectángulo áureo en subrectángulos áureos y realizaremos algunos cálculos de interés. Utilizaremos las siguientes propiedades del número fi

$$\varphi^2 = \varphi + 1 \iff \varphi = \frac{1}{\varphi} + 1 \iff \varphi - 1 = \frac{1}{\varphi}$$

**CASO I**. Recordemos que si $R$ es un rectángulo áureo de lados $a$ y $b$ entonces, $\frac{a}{b} = \varphi$ o equivalentemente $a = b\varphi$. Nótese además que:

$$a - b = b\varphi - b = \frac{b}{\varphi}$$

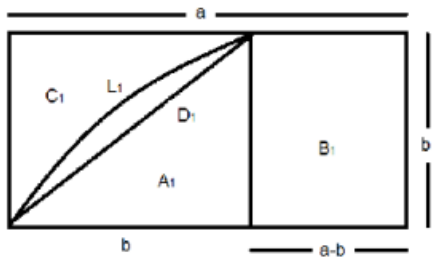

,

$$L_1 = \frac{\pi}{2}b$$

$$A_1 = \frac{\pi}{4}b^2$$

$$B_1 = (a-b)b = \frac{b^2}{\varphi}$$

$$C_1 = b^2 - \pi\frac{b^2}{4} = b^2\left(1 - \frac{\pi}{4}\right)$$

$$D_1 = \sqrt{2}b$$

Ahora consideremos el rectángulo:

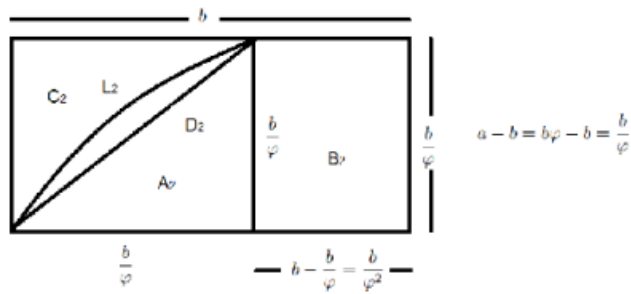

Para el cual tenemos,

$$L_2 = \frac{\pi}{2}\frac{b}{\varphi}$$

$$A_2 = \frac{\pi}{4}\frac{b^2}{\varphi^2}$$

$$B_2 = \frac{b}{\varphi}\frac{b}{\varphi^2} = \frac{b^2}{\varphi^3}$$

$$C_2 = \frac{b^2}{\varphi^2}\left(1 - \frac{\pi}{4}\right)$$

$$D_2 = \sqrt{2}\frac{b}{\varphi}$$

Siguiendo el proceso con el rectángulo:

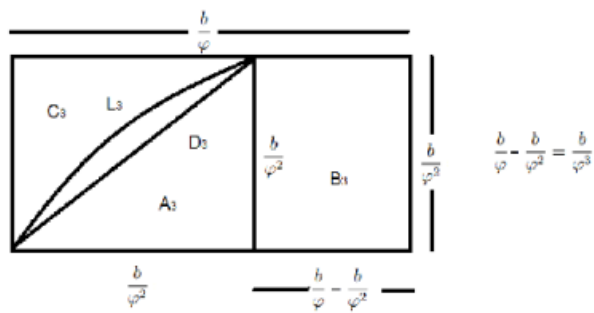

encontramos,

$$L_3 = \frac{\pi}{2}\frac{b}{\varphi^2}$$

$$A_3 = \frac{\pi}{4}\frac{b^2}{\varphi^4}$$

$$B_3 = \frac{b}{\varphi^3}\frac{b}{\varphi^2} = \frac{b^2}{\varphi^5}$$

$$C_3 = \frac{b^2}{\varphi^4}\left(1 - \frac{\pi}{4}\right)$$

$$D_3 = \sqrt{2}\frac{b}{\varphi^2}$$

Después de $k$ pasos encontramos el k-esimo rectángulo

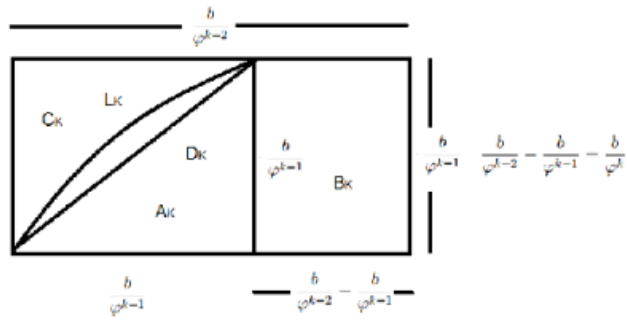

donde

$$L_k = \frac{\pi}{2}\frac{b}{\varphi^{k-1}}$$

$$A_k = \frac{\pi}{4}\frac{b^2}{\varphi^{2k-2}}$$

$$B_k = \frac{b}{\varphi^k}\frac{b}{\varphi^{k-1}} = \frac{b^2}{\varphi^{2k-1}}$$

$$C_k = \frac{b^2}{\varphi^{2k-2}}\left(1-\frac{\pi}{4}\right)$$

$$D_k = \sqrt{2}\frac{b}{\varphi^{k-1}}$$

Con estos datos podemos calcular los siguientes resultados:

$$L = \sum_{k=1}^{\infty} L_k = \frac{\pi b}{2} \sum_{k=1}^{\infty} \left(\frac{1}{\varphi}\right)^{k-1} = \frac{\pi b \varphi^2}{2} = L_1 \varphi^2$$

$$A = \sum_{k=1}^{\infty} A_k = \frac{\pi b^2}{4} \sum_{k=1}^{\infty} \left(\frac{1}{\varphi^2}\right)^{k-1} = \frac{\pi b^2 \varphi}{4} = A_1 \varphi$$

$$B = \sum_{k=1}^{\infty} B_k = b^2 \sum_{k=1}^{\infty} \frac{1}{\varphi^{2k-1}} = \frac{b^2}{\varphi} \sum_{k=1}^{\infty} \left(\frac{1}{\varphi^2}\right)^{k-1} = b^2 = B_1 \varphi$$

$$C = \sum_{k=1}^{\infty} C_k = b^2\left(1 - \frac{\pi}{4}\right) \sum_{k=1}^{\infty} \frac{1}{\varphi^{2k-2}} = b^2\left(1 - \frac{\pi}{4}\right)\varphi = C_1 \varphi$$

$$D = \sum_{k=1}^{\infty} D_k = \sqrt{2} b \sum_{k=1}^{\infty} \frac{1}{\varphi^{k-1}} = \sqrt{2} b \varphi^2 = D_1 \varphi^2$$

**CASO II**. Como caso particular tomemos $b = 1$, es decir, $a = \varphi$

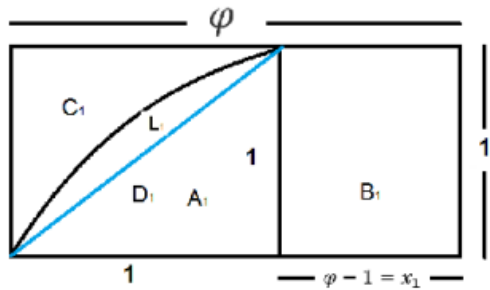

En este caso vemos que

$$L_1 = \frac{\pi}{2}, \quad A_1 = \frac{\pi}{4}, \quad B_1 = x_1 = \varphi - 1 = \frac{1}{\varphi}$$

Consideremos ahora el rectángulo:

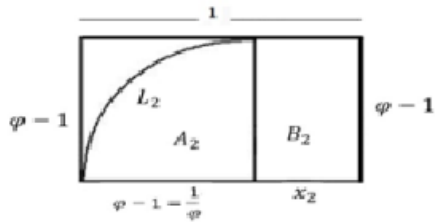

Notemos que

$$x_2 = 1 - \varphi + 1 = 1 - \frac{1}{\varphi} = \frac{\varphi - 1}{\varphi} = \frac{1}{\varphi^2}$$

También tenemos que

$$L_2 = \frac{\pi}{2}\frac{1}{\varphi}, \quad A_2 = \frac{\pi}{4\varphi^2}, \quad B_2 = (2 - \varphi)(\varphi - 1) = \frac{1}{\varphi^3}$$

Aquí tenemos un resultado importante, el área del círculo de radio $\frac{1}{\varphi}$ viene dada por

$$A = \frac{\pi}{\varphi^2} = 1{,}19998\ldots \approx 1{,}2 = \frac{6}{5}$$

De donde concluimos que

$$\varphi^2 \approx \frac{5\pi}{6} \Leftrightarrow \varphi \approx \sqrt{\frac{5\pi}{6}} \Leftrightarrow \varphi \approx \frac{5\pi}{6} - 1$$

**CASO III**. Para generalizar el caso anterior, consideremos el rectángulo de lados $a$ y $b$ con de tal manera que $a > b$ y $\frac{a}{b} = m$. Queremos saber que condición debe cumplir $m$ para que haya una división del rectángulo dado en infinitos rectángulos con la misma proporción. Notemos la sucesión de datos especiales que obtenemos,

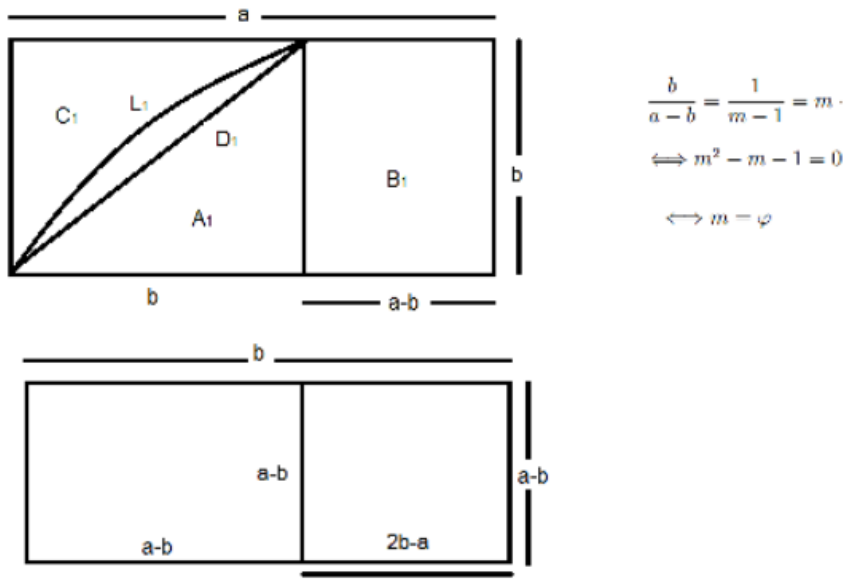

Es decir, si $\frac{a}{b} = m$ entonces en el primer subrectángulo tenemos,

$$\frac{b}{a-b} = \frac{1}{m-1} = m \iff m^2 - m - 1 = 0 \iff m = \varphi$$

Es decir, los rectángulos áureos son los únicos que permiten que sean divididos en subrectángulos áureos.

Consideremos el rectángulo de lados $a, b$ de tal manera que $a = bm, m > 1$. En el siguiente proceso podemos notar que en los subrectángulos obtenidos aparece el factor $b$.

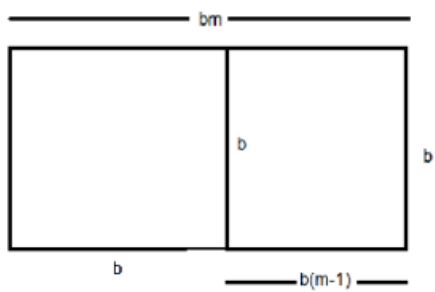
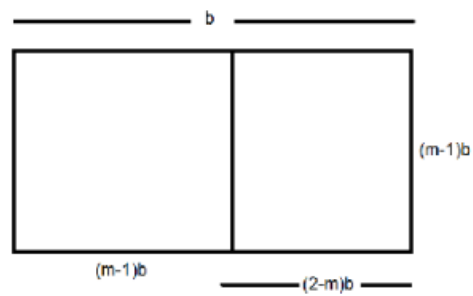
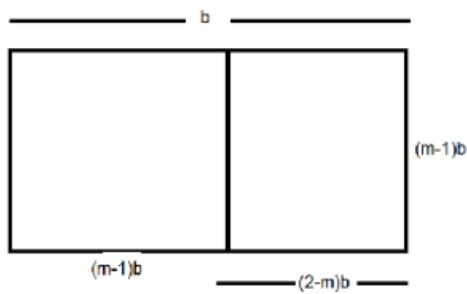
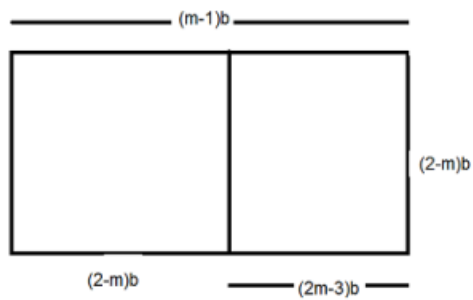

Debido a esto y por simplicidad tomaremos $b = 1, a = m$. Veremos que sucede cuando hagamos una subdivisión de k-pasos. Recuérdese que para que el proceso sea indefinido, $m = \varphi$ .

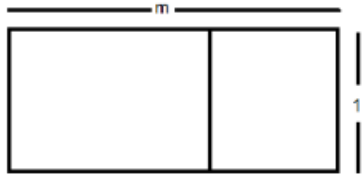

Vemos que

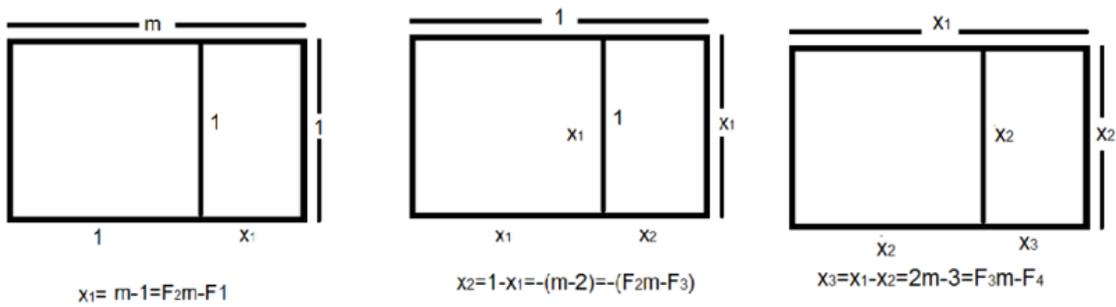

$x_1 = m-1 = F_2 m - F_1$

$x_2 = 1 - x_1 = -(m-2) = -(F_2 m - F_3)$

$x_3 = x_1 - x_2 = 2m - 3 = F_3 m - F_4$

Después de k-pasos encontramos,

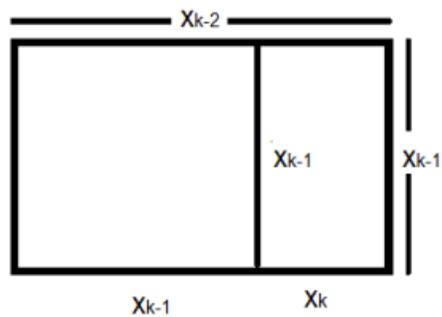

$x_k = (-1)^{k-1} (F_k m - F_{k+1})$ con $k = 1, 2, 3 \ldots$

$x_{-1} = m, \quad x_0 = 1$

$x_1 = m - 1$

Es decir,
$$x_k = (-1)^{k-1}[F_k m - F_{k+1}], \quad k = 1, 2, \ldots$$

donde $x_0 = 1, x_1 = m - 1$. Con esta fórmula podemos calcular la longitud de la espiral, de la diagonal, y las áreas mencionadas en el proceso.

Tengamos en cuenta que: para $k = 1, 2, \ldots$

$$L_k = \frac{\pi}{2} x_{k-1}.$$

$$A_k = \frac{\pi}{4} x^2\_(k-1).$$

$$B_k = x_{k-1} \cdot x\_k$$

$$C_k = x^2\_(k-1)\left(1 - \frac{\pi}{4}\right)$$

$$D_k = \sqrt{2} x_{k-1}$$

Calculemos las siguientes sumas:

$$\begin{aligned}
L &= \frac{\pi}{2} \sum_{k=1}^{n} x_{k-1} \\
&= \frac{\pi}{2}\left[x_0 + \sum_{k=2}^{n} x_{k-1}\right] \\
&= \frac{\pi}{2}\left[1 + \sum_{k=1}^{n-1} x_k\right] \\
&= \frac{\pi}{2}\left[1 + \sum_{k=1}^{n-1} (-1)^{k-1}[F_k m + F_{k+1}]\right] \\
&= \frac{\pi}{2}\left[1 + m\sum_{k=1}^{n-1}(-1)^{k-1} F_k + \sum_{k=1}^{n-1}(-1)^{k-1} F_{k+1}\right]
\end{aligned}$$

Ahora tenemos,

$$\sum_{k=1}^{n-1}(-1)^{k-1}F_k = F_1 - F_2 + \cdots + (-1)^{n-2}F_{n-1}$$

Pero sabemos que

$$F_1 - F_2 + \cdots + (-1)^{n-1}F_n = (-1)^{n-1}F_{n-1} + 1$$

Por lo que

$$F_1 - F_2 + \cdots + (-1)^{n-2}F_{n-1} = (-1)^{n-2}F_{n-2} + 1$$

También

$$\sum_{k=1}^{n-1}(-1)^{k-1}F_{k+1} = -F_2 + F_3 + \cdots + (-1)^{n-1}F_n$$

Pero como

$$F_1 - F_2 + \cdots + (-1)^{n-1}F_n = (-1)^{n-1}F_{n-1} + 1$$

Osea

$$-F_2 + \cdots + (-1)^{n-1}F_n = (-1)^{n-1}F_{n-1}$$

Reemplazando encontramos,

$$L = \frac{\pi}{2}[1 + m + (-1)^n[mF_{n-2} - F_{n-1}]]$$

Un problema sencillo pero interesante es mostrar que si $m = \varphi$ entonces

$$[1 + m + (-1)^n[mF_{n-2} - F_{n-1}]] \to \varphi^2$$

cuando $n \to \infty$

2. Calculemos ahora

$$A = \sum_{k=1}^{n} A_k = \frac{\pi}{4} \sum_{k=1}^{n} x_{k-1}^2 = \frac{\pi}{4}\left[1 + \sum_{k=1}^{n-1} x_k^2\right]$$

Pero $x_k^2 = F_k^2 m^2 - 2mF_k \cdot F_{k+1} + F_{k+1}^2$

Sumando tenemos,

$$A = \frac{\pi}{4}\left[1 + m^2 \sum_{k=1}^{n-1} F_k^2 - 2m \sum_{k=1}^{n-1} F_k F_{k+1} + \sum_{k=1}^{n-1} F_{k+1}^2\right]$$

Como $\sum_{k=1}^{n} F_k^2 = F_n F_{n+1}$ entonces

$$\sum_{k=1}^{n-1} F_k^2 = F_{n-1}F_n, \quad \sum_{k=1}^{n-1} F_{k+1}^2 = F_{n+1}F_n - 1$$

Es decir,

$$A = \frac{\pi}{4}[1 + m^2 F_{n-1}F_n - 2m \sum_{k=1}^{n-1} F_k F_{k+1} + F_{n+1}F_n - 1]$$

Como

$$\sum_{k=1}^{2n-1} F_k F_{k+1} = F_{2n}^2 \quad y \quad \sum_{k=1}^{2n} F_k F_{k+1} = F_{2n+1}^2 - 1$$

Debemos de considerar dos casos: $n = 2n$ y $n = 2n + 1$.

Para el primer caso:

$$A_{2n} = \frac{\pi}{4}[m^2 F_{2n-1} F_{2n} - 2m F_{2n}^2 + F_{2n} F_{2n+1}]$$

Si $n = 2n + 1$ entonces,

$$A_{2n+1} = \frac{\pi}{4}[m^2 F_{2n} F_{2n+1} - 2m[F_{2n+1}^2 - 1] + F_{2n+1} F_{2n+2}]$$

Nótese que si hacemos $m = \frac{F_{2n}}{F_{2n-1}}$ tenemos que

$$\begin{aligned} A_{2n} &= \frac{\pi}{4}\left[\frac{F_{2n}^2}{F_{2n-1}^2} F_{2n-1} F_{2n} - 2\frac{F_{2n}}{F_{2n-1}} F_{2n}^2 + F_{2n} F_{2n+1}\right] \\ &= \frac{\pi}{4}\left[\frac{F_{2n}^3}{F_{2n-1}^2} - 2\frac{F_{2n}^3}{F_{2n-1}} + F_{2n} F_{2n+1}\right] \\ &= \frac{\pi}{4}\left[F_{2n} F_{2n+1} - \frac{F_{2n}^3}{F_{2n-1}}\right] \\ &= \frac{\pi}{4}\left[\frac{F_{2n} F_{2n+1} F_{2n-1} - F_{2n}^3}{F_{2n-1}}\right] \end{aligned}$$

Se sabe que $F_{n-1} F_{n+1} = F_n^2 + (-1)^n$, por lo que

$$A_{2n} = \frac{\pi}{4} \frac{F_{2n}(F_{2n}^2 + 1)) - F_{2n}^3}{F_{2n-1}} = \frac{\pi}{4} \frac{F_{2n}}{F_{2n-1}}$$

Vemos que si $n \to \infty$ entonces, $A = A_{2n} = \frac{\pi\varphi}{4}$

Para $n = 2n+1$ hacemos $m = \frac{F_{2n+1}}{F_{2n}}$. Así tenemos,

$$A_{2n+1} = \frac{\pi}{4}\left[\frac{F_{2n+1}^3}{F_{2n}} - 2\frac{F_{2n+1}^3}{F_{2n}} + 2\frac{F_{2n+1}}{F_{2n}} + F_{2n+1}F_{2n+2}\right]$$

$$= \frac{\pi}{4}\left[\frac{2F_{2n+1} + F_{2n}F_{2n+1}F_{2n+1} - 2F_{2n+1}^3}{F_{2n}}\right]$$

Pero $F_{(2n)} F_{2n+2} = F_{2n+1}^2 - 1$, luego

$$A_{2n+1} = \frac{\pi}{4}\left[\frac{2F_{2n+1} + F_{2n+1}^3 - F_{2n+1} - 2F_{2n+1}^3}{F_{2n}}\right] = \frac{\pi}{4}\frac{F_{2n+1}}{F_{2n}}$$

Una vez más, si $n \to \infty$ entonces $A = A_{2n+1} = \frac{\pi\varphi}{4}$.

3. Calculemos la suma para $B_k$. Para tal efecto es conveniente tener en cuenta los siguientes resultados que enumeramos a continuación.

a)
$$\sum_{k=1}^{n} F_{n+1}^2 = F_n F_{n+1}$$

Por lo que

$$\sum_{k=1}^{n-1} F_k^2 = F_{n-1}F_n$$

$$\sum_{k=1}^{n-1} F_{k+1}^2 = \sum_{k=2}^{n} F_n F_{n+1} - 1$$

b)
$$\sum_{k=1}^{2n-1} F_k F_{k+1} = F_{2n}^2, \quad \sum_{k=1}^{2n} F_k F_{k+1} = F_{2n+1}^2 - 1$$

En particular

$$\tau(n) = \sum_{k=1}^{n-1} F_k F_{k+1} = \begin{cases} F_{2n}^2 (n = 2n) \\ F_{2n+1}^2 - 1 (n = 2n+1) \end{cases}$$

c)
$$\sum_{k=1}^{n-1} F_k F_{k+2} = \sum_{k=2}^{n} F_{k-1} F_{k+1} = \sum_{k=2}^{n} \left( F_k^2 + (-1)^k \right)$$

Ya que $F_{k-1} F_{k+1} = F_k^2 + (-1)^n$. Por lo que

$$\sum_{k=1}^{n-1} F_k F_{k+2} = F_n F_{n+1} - 1 + \sum_{k=2}^{n} (-1)^k = F_n F_{n+1} + \sum_{k=1}^{n} (-1)^k$$

d)
$$\sum_{k=1}^{n-1} F_{k+1} F_{k+2} = \sum_{k=2}^{n} F_k F_{k+1} = \sum_{k=1}^{n-1} F_k F_{k+1} + F_n F_{n+1} - 1$$

Recordemos como en los casos anteriores que

$$x_0 = 1, \quad x_1 = m - 1, x_k = (-1)^{k-1}(F_k m - F_{k+1})$$

Sabemos que $\quad B = \sum_{k=1}^{n} B_k, \quad B_k = x_{k-1} x_k$. Es decir,

$$B = \sum_{k=1}^{n} x_{k-1} x_k = x_0 x_1 + \sum_{k=1}^{n-1n} x_{k-1} x_k$$

Pero
$$x_{k-1}x_k = (-1)^{2k-1}[m^2 F_k F_{k+1} - m(F_k F_{k+2} + F_{k+1}^2) + F_{k+1}F_{k+2}]$$

Así que

$$B = m - 1 - m^2 \sum_{k=1}^{n-1} F_k F_{k+1} + m \sum_{k=1}^{n-1}(F_k F_{k+2} + F_{k+1}^2) - \sum_{k=1}^{n-1} F_{k+1} F_{k+2}$$

$$= m - \tau(n)(m^2 + 1) + m(2F_n F_{n+1} - 1 + \sum_{k=1}^{n}(-1)^k) - F_n F_{n+1}$$

donde

$$\tau(n) = \sum_{k=1}^{n-1} F_k F_{k+1}, \quad \tau(2n) = F_{2n}^2, \quad \tau(2n+1) = F_{2n+1}^2 - 1$$

Calculemos ahora para $n = 2n$ y para $n = 2n + 1$

a)
$$B_{2n} = m - F_{2n}^2(m^2 + 1) + m(2F_{2n}F_{2n+1} - 1) - F_{2n}F_{2n+1}$$
$$= -F_{2n}^2(m^2 + 1) + 2mF_{2n}F_{2n+1} - F_{2n}F_{2n+1}$$

ya que $\sum_{k=1}^{2n}(-1)^k = 0$.

Nótese en particular que si $m = \frac{F_{2n+1}}{F_{2n}}$ entonces,

$$B_{2n} = -F_{2n+1}^2 - F_{2n}^2 + 2F_{2n+1}^2 - F_{2n}F_{2n+1}$$
$$= F_{2n+1}^2 - F_{2n}(F_{2n} + F_{2n+1}) = F_{2n+1}^2 - (F_{2n+1}^2 - 1) = 1$$

Vemos que si $n \to \infty$ entonces $B = B_{2n} = 1$.

b)

$$B_{2n+1} = m - (F_{2n+1}^2 - 1)(m^2 + 1) + m(2F_{2n+1}F_{2n+2} - 2) - F_{2n+1}F_{2n+2}$$

ya que $\sum_{k=1}^{2n+1}(-1)^k = -1$. Puede mostrarse que si hacemos $m = F_{2n+1}/F_{2n+2}$ y hacemos que $n \to \infty$ entonces, $B = B_{2n+1} = 1$

4. Hagamos los cálculos finales.

Sea

$$C = \sum_{k=1}^n C_k = (1 - \frac{\pi}{4}) \sum_{k=1}^n x_{k-1}^2 = (1 - \frac{\pi}{4})\frac{4}{\pi}A = (\frac{4}{\pi} - 1)A$$

ya que $A = \frac{\pi}{4}\sum_{k=1}^n x_{k-1}^2$. Como

$$A = \frac{\pi}{4}[m^2 F_{n-1}F_n - 2m\sum_{k=1}^{n-1} F_k F_{k+1} + F_{n+1}F_n]$$

se consideran los casos $n = 2n$ y $n = 2n + 1$. Notar lo que sucede cuando $n \to \infty$.

Sea

$$D = \sum_{k=1}^n D_k = \sqrt{2}\sum_{k=1}^n x_{k-1} = \frac{2\sqrt{2}}{\pi}L$$

donde $L = \frac{\pi}{2}[1 + m + (-1)^n[mF_{n-2} + F_{n-1}]]$. De igual manera, interesa calcular $n \to \infty$

## CONCLUSIONES FINALES

Como resumen de este trabajo, tenemos las siguientes:

1. Si las medidas de los de un rectángulo son $a, b$ de tal manera que $\frac{a}{b} = m$ con 1<m<2 encontramos que para que el rectángulo se pueda subdivir en infinitos subrectángulos, tiene que ser $m = \varphi$. El rectángulo se puede dividir en $k$ subrectángulos en tanto $(-1)^{k-1}(F_k m - F_{k+1})$ sea positiva y $F_k$ son números de Fibonacci. Si $m = \frac{F_k}{F_{k-1}}$ para $k \geq 3$ el rectángulo se asemeja a un rectángulo áureo. Podemos llamar a $k$ la dimensión de áureo en un rectángulo.

2. Si un rectángulo tiene lados de medidas $a = bm$ y $b$ siendo $m > 1$, las medidas de todos los lados de los subrectángulo dependen de $m$ y $b$, por eso no se pierde generalidad al considerar los rectángulos de lados $m > 1$ y $b = 1$.

3. Si un rectángulo tiene lados de medidas $a = b\varphi$ y $b$ las siguientes son las medidas de las regiones de interes:

$$L = \sum_{k=1}^{\infty} L_k = \frac{\pi b \varphi^2}{2} = L_1 \varphi^2$$

$$A = \sum_{k=1}^{\infty} A_k = \frac{\pi b^2 \varphi}{4} = A_1 \varphi$$

$$B = \sum_{k=1}^{\infty} B_k = b^2 = B_1 \varphi$$

$$C = \sum_{k=1}^{\infty} C_k = b^2 \left(1 - \frac{\pi}{4}\right) \varphi = C_1 \varphi$$

$$D = \sum_{k=1}^{\infty} D_k = \sqrt{2b}\,\varphi^2 = D_1\varphi^2$$

4. Si en un rectángulo los lados satisfacen la relación $a = bm$, $m > 1$ se obtienen los siguientes resultados

$$L = \frac{\pi}{2}[1 + m + (-1)^n[mF_{n-2} - F_{n-1}]]$$

$$A = \frac{\pi}{4}\left[1 + m^2 F_{n-1}F_n - 2m\tau(n) + F_{n+1}F_n - 1\right]$$

$$B = -\tau(n)(m^2 + 1) + m(2F_n F_{n+1} - 1 + \sum_{k=1}^{n}(-1)^k) - F_n F_{n+1}$$

donde

$$\tau(n) = \sum_{k=1}^{n-1} F_k F_{k+1},\ \tau(2n) = F_{2n}^2,\ \ \tau(2n+1) = F_{2n+1}^2 - 1$$

$$C = \sum_{k=1}^{n} C_k = (1 - \frac{\pi}{4})\sum_{k=1}^{n} x_{k-1}^2 = (1 - \frac{\pi}{4})\frac{4}{\pi}A = (\frac{4}{\pi} - 1)A$$

$$D = \sum_{k=1}^{n} D_k = \sqrt{2}\sum_{k=1}^{n} x_{k-1} = \frac{2\sqrt{2}}{\pi}L$$

Eligiendo $m = \frac{F_h}{F_j}$ como números de Fibonacci apropiados en cada caso, se reducen a los obtenidos en 3) cuando $n \to \infty$.

## Referencias bibliográficas